\documentclass[11pt,a4paper]{article}

\usepackage{array}
\usepackage{multirow}

\makeatletter
\newcommand{\thickhline}{%
    \noalign {\ifnum 0=`}\fi \hrule height 1pt
    \futurelet \reserved@a \@xhline
}
\newcolumntype{"}{@{\hskip\tabcolsep\vrule width 1pt\hskip\tabcolsep}}
\makeatother

\usepackage{graphics}
\usepackage{epsfig}
\usepackage{subfigure}
\usepackage{amssymb}
\usepackage{amsthm}
\usepackage{mathrsfs}
\usepackage{hhline}
\usepackage{longtable}
\usepackage{amsmath}
\usepackage{float}
\usepackage{picinpar}%for window
\usepackage{cite}
\usepackage{enumerate}
\usepackage{xcolor}

\usepackage[mathlines]{lineno}
\modulolinenumbers[2]
\newcommand*\patchAmsMathEnvironmentForLineno[1]{%
\expandafter\let\csname old#1\expandafter\endcsname\csname #1\endcsname  \expandafter\let\csname oldend#1\expandafter\endcsname\csname end#1\endcsname  \renewenvironment{#1}%
{\linenomath\csname old#1\endcsname}%
{\csname oldend#1\endcsname\endlinenomath}}%
\newcommand*\patchBothAmsMathEnvironmentsForLineno[1]{%
\patchAmsMathEnvironmentForLineno{#1}%
\patchAmsMathEnvironmentForLineno{#1*}}%
\AtBeginDocument{%
\patchBothAmsMathEnvironmentsForLineno{equation}%
\patchBothAmsMathEnvironmentsForLineno{align}%
\patchBothAmsMathEnvironmentsForLineno{flalign}%
\patchBothAmsMathEnvironmentsForLineno{alignat}%
\patchBothAmsMathEnvironmentsForLineno{gather}%
\patchBothAmsMathEnvironmentsForLineno{multline}%
}

\newtheorem{theorem}{Theorem}[section]
\numberwithin{equation}{section}

\def\N{\mathbb{N}}
\newtheorem{defi}{Definition}[section]

\usepackage{graphicx}
\graphicspath{%
    {converted_graphics/}% inserted by PCTeX
    {/}% inserted by PCTeX
}

\def\N{\mathbb{N}}
\begin{document}
\baselineskip18truept
\normalsize
\begin{center}
{\mathversion{bold}\Large \bf On Integer Sequences Associated To Two Distinct Sums}

\bigskip
{\large  Gee-Choon Lau}\\

\medskip

\emph{Faculty of Computer \& Mathematical Sciences,}\\
\emph{Universiti Teknologi MARA (Segamat Campus),}\\
\emph{85000, Johor, Malaysia.}\\
\emph{geeclau@yahoo.com}\\

\end{center}

%{\blue Shiu}\quad {\red Lau}\quad {\violet Shiu 2}\\ {\magenta Shiu or Question} \quad {\cyan Maybe no use}

%\maketitle

\begin{abstract}
In this note, we show the existence of integer sequences of lengths at least 3 (except 7) such that for every integer in position $i\equiv 1\pmod{4}$ (respectively position $j\equiv 3\pmod{4}$), counting from left to right, the sum of the integer and the adjacent integer(s) has a constant sum $x$ (respectively $y$) with $x\ne y$. \\
% Include keywords, PACS and mathematical subject classification numbers as needed.


\end{abstract}

\section{Introduction}
For $m\ge 2$, let $[1,m]$ denote the set of integers from 1 to $m$. Let $\Pi_m$ be a permutation of integers in $[1,m]$, $m\ge 3$, such that for every integer in position $i\equiv 1\pmod{4}$ (respectively position $j\equiv 3\pmod{4}$), counting from left to right, the sum of the integer and the adjacent integer(s) has a constant sum $x$ (respectively $y$) with $x\ne y$. \\

Consider a graph $G(V,E)$ of order $p$ and size $q$. Let $f: V(G)\cup E(G) \to [1,p+q]$ be a bijective total labeling that induces a vertex labeling $w : V(G) \to \N$ such that $$w(v)=f(v) + \sum_{uv\in E(G)} f(uv)$$ and is called the weight of $v$ for each vertex $v \in V(G)$. We say $f$ is a local antimagic total labeling of $G$ (and $G$ is local antimagic total) if $w(u)\neq w(v)$ for each $uv\in E(G)$. Clearly, $w$ corresponds to a proper vertex coloring of $G$ if each vertex $v$ is assigned the color $w(v)$. For a graph $G$ that admits a local antimagic total labeling, the smallest number of distinct vertex weights induces by $f$ is called the local antimagic total chromatic number of $G$, denoted $\chi_{lat}(G)$ (see~\cite{Lau}). \\% In~\cite{Lau}, the following theorem is obtained. For completeness, the proof is given.

%\begin{theorem}\label{thm-LAT} Every graph $G$ is local antimagic total. \end{theorem}

%\begin{proof} It is obvious that all graphs $G$ of order $p\le 3$ are local antimagic total. We now assume $G$ is of order $p\ge4$. In~\cite{Haslegrave}, the author proved that every graph without isolated edges (by definition, necessarily without isolated vertices) admits a local antimagic labeling. Let $g$ be a local antimagic labeling of $K_1\vee G$. Define a total labeling $f : V(G) \cup E(G) \to [1, p+q]$ of $G$ such that $f(e) = g(e)$ for each edge $e\in E(G)$ and $f(v_i) = g(uv_i)$. Clearly, $w(v_i) = g^+(v_i)$ with $w(v_i)\neq w(v_j)$ if $v_iv_j\in E(G)$. Thus, $f$ is a local antimagic total labeling of $G$. \end{proof} By definition and Theorem~\ref{thm-LAT}, $\chi_{lat}(P_n)\ge 2$. 

Let $P_n = v_1v_2\cdots v_n$ be the path of order $n\ge 2$. A total labeling using integers in $[1,m], m=2n-1$ to $P_n$ such that $\chi_{lat}(P_n)=2$ corresponds to a required $\Pi_m$, $m\ge 3$ is odd. Moreover, a total labeling using integers in $[1,m], m=2n$ to $P_{n+1}$ with $v_{n+1}$ deleted such that $w(v_1) = w(v_3) = \cdots = w(v_{n-1})=x$ and $w(v_2) = w(v_4) = \cdots = w(v_{n})=y$ corresponds to a required $\Pi_m$, $m\ge 4$ is even. Obviously, a required $\Pi_3$ and $\Pi_4$ exist. In this note, we prove the following theorem. \\

\begin{theorem}\label{thm-Pi} For $m\ge 3$, there exist a required $\Pi_m$ except $m=7$.  \end{theorem}

\section{Proof}

In~\cite{Lau}, the author proved that $\chi_{lat}(P_n)=2$ if $n\ne 4$ is even and $\chi_{lat}(P_4)=3$. It was conjectured that $\chi_{lat}(P_n) = 2$ for all odd $n$. In~\cite{Lau+Shiu+Matteo+Ng}, the authors answered the conjecture in affirmative that answers Theorem~\ref{thm-Pi} for odd $m\ge 5$. For completeness, the proof is given.

\begin{theorem}\label{thm-lat-Pn} For $n\ge 2$, $\chi_{lat}(P_n)=2$ except that $\chi_{lat}(P_4)=3$. 
\end{theorem}

We consider $m\ge 5$. Theorem~\ref{thm-lat-Pn} implies that Theorem~\ref{thm-Pi} holds for odd $m\ge 5$ except $m=7$. We shall need some of the labeling functions used in proving Theorem~\ref{thm-lat-Pn}. For completeness, the functions for odd $m$ are given as follows in 4 cases.

%For $m=7$, a local antimagic total labeling is given by sequence $7$, $2$, $6$, $3$, $5$, $1$, $4$. %$P_4=v_1v_2v_3v_4$ admits a local antimagic total labeling with $f$ with $w(v_1) = w(v_3) = a$ and $w(v_2)=w(v_4)=b$. Clearly, $28 = 2a + f(v_2) + f(v_4) = 2b + f(v_1) + f(v_3)$. Thus, $f(v_2)\cong f(v_4)\pmod{2}$ and $f(v_1)\cong f(v_3)\pmod{2}$. 

%Assume $n\ge 6$ is even. Observe that by Theorem~\ref{thm-K1VG} and the local antimagic labeling of $W_n$ obtained in~\cite{Arumugam, LNS}, we can get a local antimagic total labeling of $C_n$ with an edge labeled $1$. Since $C_n$ is regular, we can delete this edge and reduce all other labels by 1. Consequently, we get a path $P_n$ that admits a local antimagic total labeling that induces exactly 2 distinct vertex weights. Since $\chi_{lat}(P_n)\ge \chi(P_n)\ge 2$, we have $\chi_{lat}(P_n)=2$. We now give a labeling that corresponds to a local antimagic total labeling of $P_n$ in~\cite{Arumugam, LNS}.

\noindent(1). For $m\equiv7\pmod{8}$, we consider $P_n, n\equiv0\pmod{4}$ (see~\cite{LNS}). A required labeling that labeled the vertices and edges of $P_{8}$ alternately is given by $10$, $7$, $12$, $1$, $13$, $3$, $11$, $6$, $9$, $2$, $14$, $4$, $8$, $5$, $15$. For $n=4k\ge 12$, we define $f: V(P_{n}) \cup E(P_{n}) \to [1, 8k-1]$ as follows.

\begin{enumerate}[(i)]
  \item $f(v_1) = 6k-2$; $f(v_{2k-1}) = 7k-1$; $f(v_{2k}) = 6k-1$; 
  \item $f(v_{2i+1}) = 4k-1+i$ for $i\in[1,k-2]$;
  \item $f(v_{2k-1+2i}) = 6k - 1 - 2i$ for $\begin{cases} \mbox{ even } k \mbox{ and } i\in[1,k/2];  \\ \mbox{ odd } k \mbox{ and } i\in[1,(k+1)/2];\end{cases}$
  \item $f(v_{3k-1+2i}) =5k - 3 + 2i$ for even $k$ and $i\in[1,k/2]$; 
  \item $f(v_{3k+2i}) = 5k - 4 + 2i$ for odd $k$ and $i\in[1,(k-1)/2]$;
  \item $f(v_{2i}) = 6k-1+i$ for $i\in[1,k-1]$;
  \item $f(v_{2k+2i}) = 8k - 2i$ for $\begin{cases} \mbox{ even } k \mbox{ and } i\in[1,k/2];  \\ \mbox{ odd } k \mbox{ and } i\in[1,(k-1)/2];\end{cases}$
  \item $f(v_{3k+2i}) = 7k-1+2i$ for even $k$ and $i\in[1,k/2]$;
  \item $f(v_{3k-1+2i}) = 7k-2+2i$ for odd $k$ and $i\in[1,(k+1)/2]$;
  \item $f(v_{2k-1}v_{2k}) = 2k-1$;
  \item $f(v_{2i-1}v_{2i}) = 3k+i$ for $i\in[1,k-1]$;
  \item $f(v_{2i}v_{2i+1}) = 2k-1-2i$ for $i\in[1,k-1]$;
  \item $f(v_{2k-1+2i}v_{2k+2i}) = 4i-2$ for $\begin{cases} \mbox{ even } k \mbox{ and } i\in[1,k/2];  \\ \mbox{ odd } k \mbox{ and } i\in[1,(k-1)/2];\end{cases}$
  \item $f(v_{3k-1+2i}v_{3k+2i}) = 2k-1+2i$ for even $k$ and $i\in[1,k/2]$;
  \item $f(v_{3k-2+2i}v_{3k-1+2i}) = 2k-2+2i$ for odd $k$ and $i\in[1,(k+1)/2]$;
  \item $f(v_{2k-2+2i}v_{2k-1+2i}) = 3k+2-2i$ for $\begin{cases} \mbox{ even } k \mbox{ and } i\in[1,k/2];  \\ \mbox{ odd } k \mbox{ and } i\in[1,(k+1)/2];\end{cases}$
  \item $f(v_{3k-2+2i}v_{3k-1+2i}) = 2k+4-4i$ for even $k$ and $i\in[1,k/2]$;
  \item $f(v_{3k-1+2i}v_{3k+2i}) = 2k+2-4i$ for odd $k$ and $i\in[1,(k-1)/2]$.
\end{enumerate}

It is not difficult to check that $$w(v_i) =\begin{cases} 9k-1 & \mbox{ for odd } i,\\ 11k-2 & \mbox{ for even } i\end{cases}$$ as required.\\ %Thus, $\chi_{lat}(P_{4k}) = 2$. %Moreover, $\cup_{k\ge 2}\{ w(v_2) - w(v_1)\} = \{3,5,7,\ldots\}$.

\begin{figure}[h] % float placement: (h)ere, page (t)op, page (b)ottom, other (p)age
  \centering
  % file name: C:/Users/geecl/Dropbox/HKNg/Local antimagic chromatic number/Chi-lat/pi23.eps
  \includegraphics[bb=0 0 201 44,width=2.56in,height=0.56in,keepaspectratio]{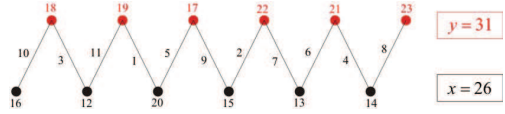}   
  \caption{$\Pi_{23}$ with $x=26$ and $y=31$}
  \label{fig:pi23}
\end{figure}

\begin{figure}[h] % float placement: (h)ere, page (t)op, page (b)ottom, other (p)age
  \centering
  % file name: C:/Users/geecl/Dropbox/HKNg/Local antimagic chromatic number/Chi-lat/pi31.eps
  \includegraphics[bb=0 0 274 45,width=3.43in,height=0.56in,keepaspectratio]{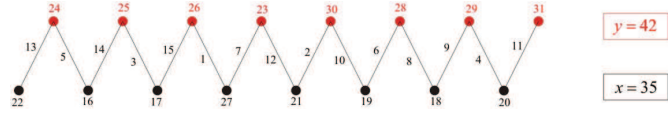}
  \caption{$\Pi_{31}$ with $x=35$ and $y=42$}
  \label{fig:pi31}
\end{figure}

The labeling sequence for $\Pi_{23}$ is $16$, $10$, $18$, $3$, $12$, $11$, $19$, $1$, $20$, $5$, $17$, $9$, $15$, $2$, $22$, $7$, $13$, $6$, $21$, $4$, $14$, $8$, $23$ with $2$ distinct vertex weights that correspond to $x=26$ and $y=31$. The labeling sequence for $\Pi_{31}$ is $22$, $13$, $24$, $5$, $16$, $14$, $25$, $3$, $17$, $15$, $26$, $1$, $27$, $7$, $23$, $12$, $21$, $2$, $30$, $10$, $19$, $6$, $28$, $8$, $18$, $9$, $29$, $4$, $20$, $11$, $31$ with $2$ distinct vertex weights that correspond to $x=35$ and $y=42$. (See Figures~\ref{fig:pi23} and~\ref{fig:pi31}.)

\noindent(2). For $m\equiv3\pmod{8}$, we consider $P_n, n\equiv2\pmod{4}$. For $n=4k+2\ge 6$, we define $f: V(P_n) \cup E(P_n) \to [1, 8k+3]$ as follows.

\begin{enumerate}[(i)]
  \item $f(v_{4i+2}) = 5k+1-i$ for $i\in [0,k-1]$;
  \item $f(v_{4i+4}) = 6k+1-i$ for $i\in [0,k-1]$;
  \item $f(v_{4k+2}) = 6k+2$;
  \item $f(v_{2i}v_{2i+1}) = i$ for $i\in [1,2k]$;
  \item $f(v_{4i+1}v_{4i+2}) = 4k+1-i$ for $i\in [0,k]$;
  \item $f(v_{4i+3}v_{4i+4}) = 3k-i$ for $i\in [0,k-1]$.
\end{enumerate}

It is not difficult to check that $$w(v_i) =\begin{cases} 11k+4 & \mbox{ for odd } i,\\ 9k+3 & \mbox{ for even } i\end{cases}$$  as required. \\

\begin{figure}[h] % float placement: (h)ere, page (t)op, page (b)ottom, other (p)age
  \centering
  % file name: C:/Users/geecl/Dropbox/HKNg/Local antimagic chromatic number/Chi-lat/pi27.eps
  \includegraphics[bb=0 0 245 45,width=3.07in,height=0.56in,keepaspectratio]{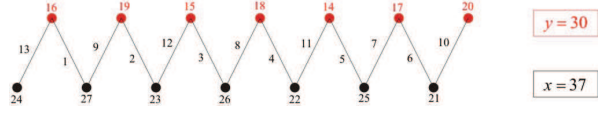}
  \caption{$\Pi_{27}$ with $x=37$ and $y=30$}
  \label{fig:pi27}
\end{figure}

% The following labeling sequences of $P_n$ corresponds to a required $\Pi_{2n-1}$.\\
The labeling sequence for $\Pi_{27}$ is $24$, $13$, $16$, $1$, $27$, $9$, $19$, $2$, $23$, $12$, $15$, $3$, $26$, $8$, $18$, $4$, $22$, $11$, $14$, $5$, $25$, $7$, $17$, $6$, $21$, $10$, $20$ with $2$ distinct vertex weights that correspond to $x=37$ and $y=30$. (See Figure~\ref{fig:pi27}.)

\noindent(3). For $m\equiv1\pmod{8}$, we consider $P_n, n\equiv1,3\pmod{4}$ (see~\cite{Arumugam}). Consider $n=4k+1\ge 5$. For $n=5$, a required labeling sequence that labeled the vertices and edges of $P_5$ alternately is $6$, $4$, $7$, $3$, $2$, $5$, $8$, $1$, $9$  with distinct vertex weights 10 and 14. By computer search, we are able to obtain all the 12 different labelings. For $k\ge 2$, we define $f : V(P_n) \cup E(P_n) \to [1, 8k+1]$ as follows.

\begin{enumerate}[(i)]
  \item $f(v_{1}) = 8k$; $f(v_{4k-1}) = 6k$; $f(v_{4k+1}) = 8k+1$; 
  \item $f(v_{4i+1}) = 7k-3+i$ for $i\in [1,k-1]$;
  \item $f(v_{4i-1}) = 3k+i$ for $i\in [1,k-1]$;
  \item $f(v_{4i+2}) = 7k+i$ for $i\in [0,k-1]$;
  \item $f(v_{4i+4}) = 4k+1+i$ for $i\in [0,k-1]$;
  \item $f(v_{4k}v_{4k+1}) = 2k$; $f(v_{4k-1}v_{4k}) = 4k$; 
  \item $f(v_{2i}v_{2i+1}) = 2k-i$ for $i\in [1,2k-1]$;
  \item $f(v_{4i+1}v_{4i+2}) = 2k+1+i$ for $i\in [0,k-1]$;
  \item $f(v_{4i+3}v_{4i+4}) = 5k+1+i$ for $i\in [0,k-2]$.
% \item $f(v_1) = 8k+1$; $f(v_3) = 6k$; $f(v_{4k+1}) = 8k$;
% \item $f(v_{4i+1}) = 7k-i$ for $i\in [1,k-1]$; %$1\le i\le k-1$;
% \item $f(v_{4i+3}) = 4k-i$ for $i\in [1,k-1]$; %$1\le i\le k-1$;
% \item $f(v_{4i+2}) = 5k-i$ for $i\in [0,k-1]$; %$0\le i\le k-1$;
% \item $f(v_{4i+4}) = 8k-1-i$ for $i\in [0,k-1]$; %$0\le i\le k-1$;
% \item $f(v_1v_2) = 2k$; $f(v_2v_3) = 4k$;
% \item $f(v_{2i+1}v_{2i+2}) = i$ for $i\in [1,2k-1]$; %$1\le i\le 2k-1$;
% \item $f(v_{4i}v_{4i+1}) = 3k+1-i$ for $i\in [1,k]$; %$1\le i\le k$;
% \item $f(v_{4i+2}v_{4i+3}) = 6k-i$ for $i\in [1,k-1]$; %$1\le i\le k-1$.
\end{enumerate}

It is not difficult to check that $$w(v_i) =\begin{cases} 10k+1 & \mbox{ for odd } i,\\ 11k & \mbox{ for even } i\end{cases}$$ as required.\\ %Moreover, $\cup_{k\ge 2}\{ w(v_2) - w(v_1)\} = \Z^+$.

\begin{figure}[h] % float placement: (h)ere, page (t)op, page (b)ottom, other (p)age
  \centering
  % file name: C:/Users/geecl/Dropbox/HKNg/Local antimagic chromatic number/Chi-lat/pi25.eps
  \includegraphics[bb=0 0 230 45,width=2.88in,height=0.56in,keepaspectratio]{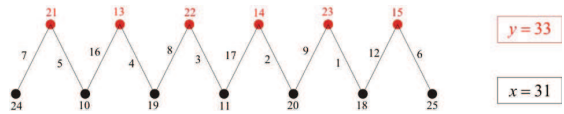}
  \caption{$\Pi_{25}$ with $x=31$ and $y=33$}
  \label{fig:pi25}
\end{figure}

The labeling sequence for $Pi_{25}$ is $24$, $7$, $21$, $5$, $10$, $16$, $13$, $4$, $19$, $8$, $22$, $3$, $11$, $17$, $14$, $2$, $20$, $9$, $23$, $1$, $18$, $12$, $15$, $6$, $25$ with $2$ distinct vertex weights that correspond to $x=31$ and $y=33$. (See Figure~\ref{fig:pi25}.)

\noindent(4). For $m\equiv5\pmod{8}$, we consider $P_n, n\equiv3\pmod{4}$. Consider $n=4k+3\ge 3$. A required labeling sequence for $n=3,7,11$ are (i) $1$, $5$, $3$, $4$, $2$; (ii) $5$, $13$, $1$, $9$, $7$, $2$, $10$, $11$, $4$, $3$, $8$, $12$, $6$; and (iii) $8$, $21$, $1$, $14$, $11$, $4$, $15$, $17$, $10$, $2$, $16$, $18$, $6$, $5$, $12$, $19$, $7$, $3$, $12$, $20$, $9$. The corresponding distinct vertex weights are 6 and 12; 18 and 23; and 29 and 36 respectively. 

For $k\ge 3$, we define $f: V(P_n)\cup E(P_n)\to [1, 8k+5]$ as follows.

\begin{enumerate}[(i)]
  \item $f(v_1) = 3k+2$; $f(v_2) = 1$; $f(v_3) = 4k+3$; $f(v_{2k+3}) = 2k+2$; $f(v_{4k+3}) = 3k+3$;
  \item $f(v_{2i+3}) = 3k+3+i$ for $i\in[1,k-1]$;
  \item $f(v_{2k+2i+1}) = 2k+1+i$ for $i\in[1,k]$;
  \item $f(v_{2i+2}) = 5k+4+i$ for $i\in[1,k]$;
  \item $f(v_{2k+2i+2}) = 4k+3+i$ for $i\in[1,k]$;
  \item $f(v_1v_2) = 8k+5$; $f(v_2v_3) = 5k+4$;
  \item $f(v_{2i+1}v_{2i+2}) = 2k+2-2i$ for $i\in[1,k]$;
  \item $f(v_{2k+2i+1}v_{2k+2i+2}) = 2k+3-2i$ for $i\in[1,k]$;
  \item $f(v_{2i+2}v_{2i+3}) = 6k+4+i$ for $i\in[1,2k]$.
\end{enumerate}

It is not difficult to check that $$w(v_i) =\begin{cases} 11k+7 & \mbox{ for odd } i,\\ 13k+10 & \mbox{ for even } i\end{cases}$$  as required.

% The labeling sequences $1,5,3,4,2$; $1,9,7,3,2,5,8,6,4$ and $13,6,4,10,1,8,9,3,5,11,2,7,12$ show that $\chi_{lat}(P_n)=2$ for $n=3,5,7$ respectively. Another labeling for $P_7$ is $6,12,8,3,4,11,10,2,7,9,1,13,5$.

\begin{figure}[h] % float placement: (h)ere, page (t)op, page (b)ottom, other (p)age
  \centering
  % file name: C:/Users/geecl/Dropbox/HKNg/Local antimagic chromatic number/Chi-lat/pi29.eps
  \includegraphics[bb=0 0 259 45,width=3in,height=0.56in,keepaspectratio]{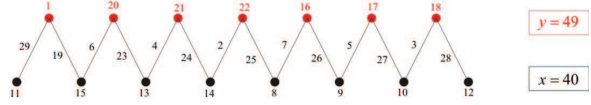}
  \caption{$\Pi_{29}$ with $x=40$ and $y=49$}
  \label{fig:pi29}
\end{figure}

The labeling sequence for $Pi_{29}$ is $11$, $29$, $1$, $19$, $15$, $6$, $20$, $23$, $13$, $4$, $21$, $24$, $14$, $2$, $22$, $25$, $8$, $7$, $16$, $26$, $9$, $5$, $17$, $27$, $10$, $3$, $18$, $28$, $12$ with $2$ distinct vertex weights that correspond to $x=40$ and $y=49$. (See Figure~\ref{fig:pi29}.)

%\medskip Suppose $F_n$ is the fan graph of order $n$ obtained from $W_n$ by deleting an outer edge of $W_n$. In~\cite[Theorem 3.7]{LSN1}, the authors showed that $3\le\chi_{la}(F_n)\le 4$ for odd $n\ge 3$. For $n=4k+1$, we have $\sum^{4k+1}_{i=1} f(u_i) = 22k^2 + 12k + 1 \ne 11k \ne 10k+1$. For $n=4k+3$, we have $\sum^{4k+3}_{i=1} f(u_i) = 16k^2 + 19k + 6 \ne 11k+7 \ne 13k+10$. By Theorem~\ref{thm-K1VG}, we have the following.

%\begin{corollary} For $k\ge 0$, $\chi_{la}(F_{2k+1}) = 3$. 
%\end{corollary}
%\begin{conjecture} For odd $n\ge 3$, $\chi_{lat}(P_n) = 2$. \end{conjecture}

%\begin{theorem}\label{thm-pi} For each $m\ge 3$, there exists a permutation sequence $\Pi_m$ as defined above except $m=7$.
%\end{theorem}
%For $n\ge 2$ and $n\ne4$, consider each local antimagic total labeling of the path $P_n$, $n = 4j+1, 4j+2, 4j+3, 4j+4$ given above. For $1\le i\le n$, let the integer assign to vertex $v_i$ be in position $2i-1$. For $1\le i\le n-1$, let the integer assign to edge $v_iv_{i+1}$ be in position $2i$. 

We have now obtained a required sequence $\Pi_m$, $m=8j+1, 8j+3, 8j+5, 8j+7$ respectively. \\ 

%
%Suppose $n\equiv0\pmod{8}$, we note that the sequence $P$ of length $n\equiv 1\pmod{8}$ begins with the largest integer. Reversing the sequence and deleting the largest integer, we get a required sequence $P$.
%
We now consider even $m\ge 6$. We note that for $m=8j+1$ and $m=8j+7$, each sequence $\Pi_m$ ends with the largest integer $m$. Deleting $m$, we get a required $\Pi_{m-1}$ of length $8j$ and $8j+6$ as required.\\

For $m=8j+2, j\ge 1$, we define a total labeling $f: V(P_{4j+2})\cup E(P_{4j+2}))\setminus\{v_{4j+2}\} \to [1,8j+2]$ as follows:

\begin{enumerate}[(i)] 
 % \item $f(v_{4j+2}) = 0$;
  \item $f(v_{4i-2}) = 3j+1-i$ for $i\in[1,j]$;
  \item $f(v_{4i}) = 4j+1-i$ for $i\in[1,j]$;
  \item $f(v_{4i-3}) = 7j+3-i$ for $i\in[1,j+1]$;
  \item $f(v_{4i-1}) = 8j+3-i$ for $i\in[1,j]$;
  \item $f(v_{2i}v_{2i+1}) = i$ for $i\in[1,2j]$;
  \item $f(v_{4i-3}v_{4i-2}) = 6j+2-i$ for $i\in[1,j+1]$;
  \item $f(v_{4i-1}v_{4i}) = 5j+1-i$ for $i\in[1,j]$.
\end{enumerate} 

It is not difficult to check that $$w(v_i) =\begin{cases} 13j+3 & \mbox{ for odd } i,\\ 9j+2 & \mbox{ for even } i.\end{cases}$$ %List the assigned integers as in Example~\ref{ex-path} and delete the ending $0$, we get a required $\Pi_{8j+2}$.\\

\begin{figure}[h] % float placement: (h)ere, page (t)op, page (b)ottom, other (p)age
  \centering
  % file name: C:/Users/geecl/Dropbox/HKNg/Local antimagic chromatic number/Chi-lat/pi26.eps
  \includegraphics[bb=0 0 245 45,width=3.07in,height=0.56in,keepaspectratio]{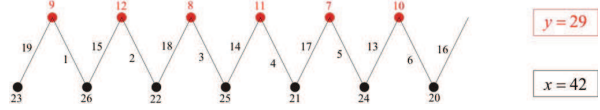}
  \caption{$\Pi_{26}$ with $x=42$ and $y=29$}
  \label{fig:pi26}
\end{figure}

The sequence $\Pi_{26}$ is given by $23$, $19$, $9$, $1$, $26$, $15$, $12$, $2$, $22$, $18$, $8$, $3$, $25$, $14$, $11$, $4$, $21$, $17$, $7$, $5$, $24$, $13$, $10$, $6$, $20$, $16$ with $x=42$ and $y=29$. (See Figure~\ref{fig:pi26}.)\\

For $m=8j+4, j\ge 1$, we define a total labeling $f: (V(P_{4j+3}) \cup E(P_{4j+3}))\setminus\{v_{4j+3}\} \to [1,8j+4]$ as follows:
\begin{enumerate}[(i)] 
 % \item $f(v_{4j+3}) = 0$;
  \item $f(v_{2i-1}) = 7j+5-i$ for $i\in[1,j+1]$;
  \item $f(v_{2j+1+2i}) = 8j+5-i$ for $i\in[1,j]$;
  \item $f(v_{2i}) = 5j+3-i$ for $i\in[1,j]$;
  \item $f(v_{2j+2i}) = 7j+1-i$ for $i\in[1,j+1]$;
  \item $f(v_{2i-1}v_{2i}) = 4j+3-i$ for $i\in[1,2j+1]$;
  \item $f(v_{2i}v_{2i+1}) = 2i$ for $i\in[1,j]$;
  \item $f(v_{2j+2i}v_{2j+2i+1}) = 2i-1$ for $i\in[1,j+1]$.
\end{enumerate}

It is not difficult to check that $$w(v_i) =\begin{cases} 11j+6 & \mbox{ for odd } i,\\ 9j+6 & \mbox{ for even } i.\end{cases}$$ %List the assigned integers as in Example~\ref{ex-path} and delete the ending $0$, we get a required $\Pi_{8j+4}$.

%\begin{example} The following labeling sequences of $P_n$ with the right end-vertex deleted correspond to a required $\Pi_{2n-2}$.\\

\begin{figure}[h] % float placement: (h)ere, page (t)op, page (b)ottom, other (p)age
  \centering
  % file name: C:/Users/geecl/Dropbox/HKNg/Local antimagic chromatic number/Chi-lat/pi20.eps
  \includegraphics[bb=0 0 201 45,width=2.52in,height=0.56in,keepaspectratio]{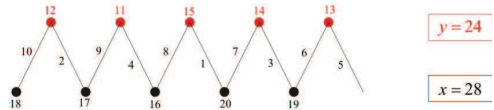}
  \caption{$\Pi_{20}$ with $x=28$ and $y=24$}
  \label{fig:pi20}
\end{figure}

The sequence $\Pi_{20}$ is given by $18$, $10$, $12$, $2$, $17$, $9$, $11$, $4$, $16$, $8$, $15$, $1$, $20$, $7$, $14$, $3$, $19$, $6$, $13$, $5$ with $x=28$ and $y=24$. (See Figure~\ref{fig:pi20}.)
%\end{example}

%{\bf NOTE:} In the permutation sequences in Theorem~\ref{thm-Pi}, we have $$|x - y| = \begin{cases}j-1\ge 1 &\mbox{ for } m = 8j, 8j+1; \\ 2j \ge 2 &\mbox{ for } m = 8j+4; \\ 2j+1\ge 3 &\mbox{ for } m = 8j+3, 8j+6, 8j+7; \\ 2j+3\ge 9 &\mbox{ for } m=8j+5; \\ 4j+1\ge 5 &\mbox{ for } m = 8j+2. \end{cases}$$

%\vskip-10cm

\begin{thebibliography}{99}
\bibitem{Arumugam} S. Arumugam, K. Premalatha, M. Bac\v{a} and A. Semani\v{c}ov\'{a}-Fe\v{n}ov\v{c}\'{i}kov\'{a}, Local antimagic vertex coloring of a graph, {\it Graphs and Combin.}, {\bf33} (2017), 275 -- 285.

%\bibitem{Bensmail} J. Bensmail, M. Senhaji and K. Szabo Lyngsie, On a combination of the 1-2-3 Conjecture and the Antimagic Labelling Conjecture, {\it Discrete Math. Theoret. Comput. Sc.}, {\bf19(1)} (2017) \#22.

%\bibitem{Chai} E.S. Chai, A. Das and C.K. Midha, Construction of magic rectangles of odd order, {\it Australas. J. Combin.}, {\bf 55} (2013), 131--144.

%\bibitem{Reyes} J.P. De Los Reyes, A. Das and C.K. Midha, A matrix approach to construct magic rectangles of even order, {\it Australas. J. Combin.}, {\bf 40} (2008), 293--300.

\bibitem{Haslegrave} J. Haslegrave, Proof of a local antimagic conjecture, {\it Discrete Math. Theor. Comp. Sc.}, {\bf 20(1)} (2018), \#18.

\bibitem{Lau} G.C. Lau, Every graph is local antimagic total, (2019) arXiv:1906.10332.

\bibitem{Lau+Shiu+Matteo+Ng} G.C. Lau, W.C. Shiu, M. Bodini, H.K. Ng, On local antimagic (total) chromatic numbers, (2020), in preparation.

\bibitem{LNS} G.C. Lau, H.K. Ng, and W.C. Shiu, Affirmative solutions on local antimagic chromatic number, (2019) arXiv:1805.02886.

%\bibitem{LauNgShiu-CM} Lau, G.C., Ng, H.K., Shiu, W.C.: Cartesian magicness of 3-dimensional boards. (2018) arXiv:1805.04890.

%\bibitem{LSN1} G.C. Lau, W.C. Shiu and H.K. Ng, On local antimagic chromatic number of cycle-related join graphs, {\it Discuss. Math. Graph Theory} (2018), doi:10.7151/dmgt.2177. 

\end{thebibliography}
\end{document}